# Assembly Procedure for Elementary Matrices of Train-Track-Bridge Railway System

**Laila Bouhlal[1*], Fouzia Kassou[2], Nouzha Lamdouar[3], Azeddine Bouyahyaoui[4]**

[1]Department of Civil Engineering, Mohammedia Engineering School, Mohammed V University, Rabat, Morocco
[2]Department of Civil Engineering, Ecole Hassania des Travaux Publics, Casablanca, Morocco
[3]Department of Civil Engineering, Mohammedia Engineering School, Mohammed V University, Rabat, Morocco
[4]Department of Civil Engineering, Mohammedia Engineering School, Mohammed V University, Rabat, Morocco
*Corresponding Author: lailabouhlal@emi.ac.ma

**Abstract** The aim of this article is to determine the dynamic response of the railway system composed of train, track and bridge, suggesting a method for assembling the elementary matrices to obtain differential equations of the overall system. The model studied consists of a moving part which is the vehicle modeled by a mass-spring-damper system, and a fixed part made up of the rail and bridge deck, modeled by two Bernoulli beams. The ballast layer of the rail bed is represented by springs and continuous dampers. The motion equations are established using the principle of total stationary potential energy in combination with the finite element method.

The Newmark method was used to solve the motion equations, whose coefficients depend on time. It is an explicit numerical integration method with unconditional stability and high accuracy, which does not involve any iterative procedure. The concept of the method is explained herein and applied to the railway system studied.

A computer program was developed on MATLAB software to analyze the system's dynamic responses. The assembly procedure presented, and the numerical method used to solve the equations in the time domain have been validated there an illustrated example in literature.

## 1. Introduction

Numerical simulation has become an effective tool for dynamic prediction of vehicle, track and bridge behavior in the railway field [1]. Several models have been developed in various studies, including [2-4], which presented the vehicle-rail-bridge in two subsystems namely the vehicle and the bridge-rail. References [5,6] consider the vehicle-rail-bridge system as a single coupled system, the dynamic contact forces between train and rail being internal forces.

Dynamic responses determination of a railway system involves 3 critical steps requiring high accuracy respectively:

- Modeling and elaboration of elementary differential equations,
- Assembling elementary matrices to establish equations of the overall system,
- Selecting a numerical method for solving equations, and consequently for predicting the dynamic behavior of the system studied.

Several methods have been used to determine motion equations, reference [7] reported the motion equations of the vehicle-rail-bridge system using Hamilton's principle. D'Alembert's principle was used [8] to establish set of differential equations.

The principle of the stationary value of the total potential energy combined with the finite element method has been adopted in the current study to establish motion equations of the system. The dynamic interaction model of the vehicle-rail-bridge system has evolved over time from a simple restricted model of constant force in motion to axle loads. Several studies have adopted the restricted model [9,10], which sought only to predict axle response without taking account track irregularities. The moving and suspended mass model remained simple after moving load and allowed us to study the effect of vehicle inertia on axle response [11,12].

Other researchers have modeled the train by two rigid masses connected by spring-damper suspension systems [13,14]. The moving systems provide more realistic modeling and behavior for the interaction issue. The tow DOFs (degrees of freedom) system includes a first mass which represents the vehicle body and the second, the wheel mass. This model has 2DOFs respectively the vertical displacement of the body and the wheel. The 2 DOFs model has been modified to incorporate the effect of pitching, with a 4 DOFs system [15]. The 4 DOFs are the two wheels vertical translations, translation and rotation body. Reference [16] developed a 10 DOFs model, the train being modeled by a body, two bogies and 4 wheels, connected respectively by a primary and a secondary suspension represented by spring-damper systems. This model has 10 DOFs and allows the study of the bridge dynamic response and passengers comfort [17,18]. Reference [19] presented a more complex 115 DOFs train model including vertical and transverse connections between car bodies and two suspension layers.

After modeling and determination of motion equations for the elementary system, an assembly process was developed and presented to generate differential equations of the overall system composed of vehicle, track and bridge. The engineer is then asked to solve these differential equations with several degrees of freedom. When the system is non-linear, or the applied force varies arbitrarily with time, the mathematical resolution seems more complex, or even impossible in some cases.

Widely known for its convergence and precision, the Newmark beta method was chosen to analyze the dynamic behavior of the railway system studied.

# 2. Equations of Motion for "Train-Track-Bridge" Interaction

## 2.1. The Theoretical Model

The vehicle rail and bridge are considered a complete system in current investigation. The vehicle is modeled by two rigid masses respectively, $m_l$ body mass and me wheel mass. The two masses are connected by a spring-damper system ($k_l$,$c_l$), $q_l$ and qe present vertical displacements of the vehicle which is running with a velocity v(t) and an acceleration a(t) in longitudinal direction. Therefore, total number of DOFs is 2, since vertical movement of wheel is limited by the rail, the number of independent DOFs becomes 1 noted $\{q_1\}$.

Rail and bridge are modeled by two Euler-Bernoulli beams of finite length simply supported on bridge piers. Two beams are connected by a continuous layer of damper springs (krb,crb). Based on the finite element method, rail and bridge are divided into 10 elements of length l, the damping of the rail is neglected, however a linear viscous damping $c_b$ is considered for the bridge, $m_b$ and $m_r$ represent respectively the mass per unit length for bridge and rail. when longitudinal displacement is neglected of the two beams, each beam element has 4 DOFs, vertical displacement and rotation in each extremity. Displacement vectors respectively for rail and bridge are $\{q_2,q_3,q_4,q_5\}$ and $\{q_6,q_7,q_8,q_9\}$.Therefore vehicle-rail-bridge element DOFs become 9 $\{q_1,q_2,q_3,q_4,q_5,q_6,q_7,q_8,q_9\}$.

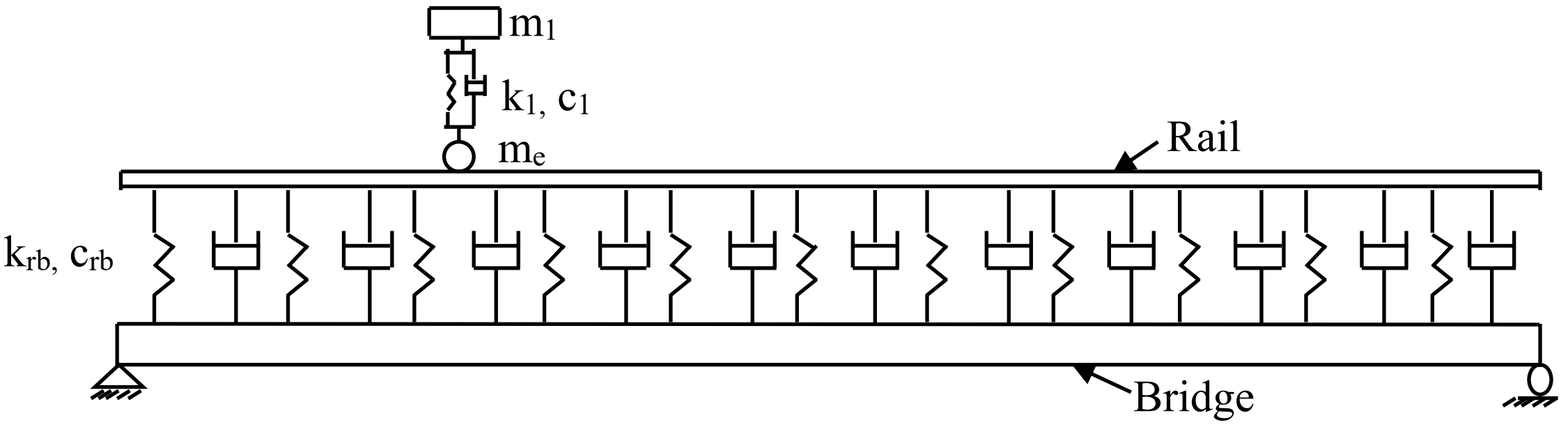

**Figure 1.** Model of Vehicle-Track-Bridge Interaction System

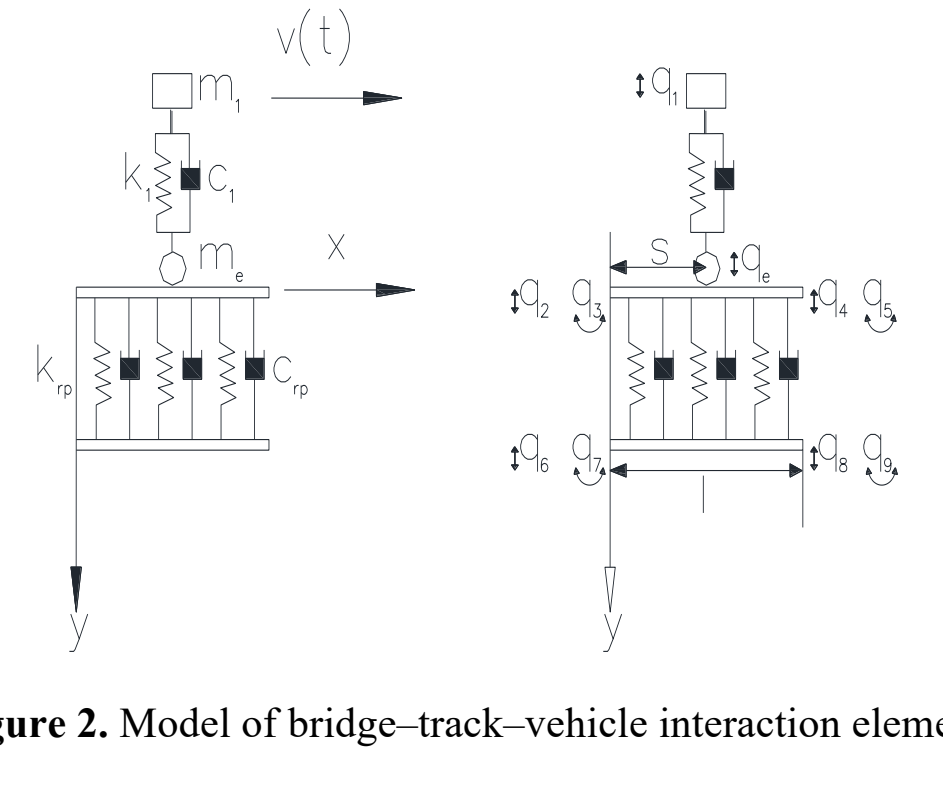

**Figure 2.** Model of bridge–track–vehicle interaction element

## 2.2. Formulation of Motion Equations

Formulation of motion equations for a dynamic system remains the key and delicate step of the dynamic analysis process. The principle of a stationary value of total potential energy is adopted herein. The potential energy of a dynamic system includes potential energy $E_{pg}$ of gravity, the elastic deformation energy $E_{pe}$, the potential energy of the damping force $E_{pa}$, the potential energy $E_{pi}$ of the inertia force [11], hence:

$$E_p=E_{pg}+E_{pe}+E_{pa}+E_{pi} \quad (1)$$

The principle of least action requires that the virtual work done when the mechanical system performs displacements is null.

$$\delta E_p = \delta W = 0 \quad (2)$$

$$\delta Ep = \delta Epg + \delta Epe + \delta Epa + \delta Epi = 0 \quad (3)$$

In a beam of length element l, the displacement w(s,t) at any point s of the beam should be described in terms of nodal displacements at beam extremities { q1,q2,q3,q4}. Interpolation functions are used. At each time t, for any point $s \in [0;l]$ we have then:

$$w(s,t) = N_A(s)q_1(t) + N_B(s)q_2(t) + N_C(s)q_3(t) + N_D(s)q_4(t) \quad (4)$$

$$w''(s,t) = N_A''(s)q_1(t) + N_B''(s)q_2(t) + N_C''(s)q_3(t) + N_D''(s)q_4(t) \quad (5)$$

$$\dot{w}(s,t) = N_A(s)\dot{q}_1(t) + N_B(s)\dot{q}_2(t) + N_C(s)\dot{q}_3(t) + N_D(s)\dot{q}_4(t) \quad (6)$$

$$\ddot{w}(s,t) = N_A(s)\ddot{q}_1(t) + N_B(s)\ddot{q}_2(t) + N_C(s)\ddot{q}_3(t) + N_D(s)\ddot{q}_4(t) \quad (7)$$

By setting $\xi = \frac{s}{l}$, they are given by:

$$N_A(\xi) = 1 - 3\xi^2 + 2\xi^3 \quad (8)$$

$$N_B(\xi) = l\xi(\xi-1)^2 \quad (9)$$

$$N_C(\xi) = \xi^2(3-2\xi) \quad (10)$$

$$N_D(\xi) = l\xi^2(\xi-1) \quad (11)$$

At a fixed time t, the wheel of mass me is located at a distance s from rail left end.

We note that: x=v×t and j is the number of the element on which the wheel acts, with : $j = E(\frac{x}{l}) + 1$

we specify that : $E(\frac{x}{l})$ the integer part of the number ,

v: the wheel speed. In this case s=x-(j-1)×l.

Taking into consideration track irregularities causing vertical deviation along the rail from its initially horizontal profile, let r(s) be the value of this deviation at the point of contact with the wheel.

The condition of permanent contact between wheel and rail results in the following:

$$q_e(s,t) = r(s) + w_r(s,t) \quad (12)$$

Where:

$$\dot{q}_e(s,t) = vr'(s) + \dot{w}_r(s,t) + vw'(s,t) \quad (13)$$

$$\ddot{q}_e(s,t) = v^2 r(s) + ar'(s) + \ddot{w}(s,t) + 2v\left[\frac{\partial^2}{\partial x \partial t} w(x,t)\right]_{x=s} + v^2 w(s,t) + aw'(s,t) \quad (14)$$

v and a are respectively the speed and acceleration of vehicle in horizontal direction.

$$N(x) = [N_A(x)\ N_B(x)\ N_c(x)\ N_D(x)]$$

Calculating the variations of energies and replacing them in "(3)", the motion equations of the elementary system vehicle-rail-bridge Fig. 2 are derived.

Therefore, the equations in matrix form can be written as follows:

$$M_e Q_e + C_e \dot{Q}_e + K_e \ddot{Q}_e = F_e \quad (15)$$

With:

$$Me = \begin{bmatrix} [m_{11}] & [0] & [0] \\ \{0\} & [m_{rr}] & [0] \\ \{0\} & [0] & [m_{pp}] \end{bmatrix}$$

$$Ke = \begin{bmatrix} [k_{11}] & [k_{1r}] & [0] \\ \{k_{r1}\} & [k_{rr}] & [k_{rp}] \\ \{0\} & [k_{pr}] & [k_{pp}] \end{bmatrix}$$

$$Ce = \begin{bmatrix} [c_{11}] & [c_{1r}] & [0] \\ \{c_{r1}\} & [c_{rr}] & [c_{rp}] \\ \{0\} & [c_{pr}] & [c_{pp}] \end{bmatrix}$$

We note that:

Vector of nodal rail displacement:

$$\{q_r\} = \begin{Bmatrix} q_2 \\ q_3 \\ q_4 \\ q_5 \end{Bmatrix}$$

Vector of nodal bridge displacement:

$$\{q_p\} = \begin{Bmatrix} q_6 \\ q_7 \\ q_8 \\ q_9 \end{Bmatrix}$$

Mass Matrix component:

$$[m_{11}] = m_1$$

$$[m_{rr}] = m_e N(s)^T N(s) + m_r \int_0^l N(x)^T N(x) dx$$

$$[m_{pp}] = m_p \int_0^l N(x)^T N(x) dx$$

Damping Matrix component:

$$[c_{11}] = c_1 \ , \ \{c_{r1}\} = \{c_{1r}\} = -c_1 N(s)$$

$$[c_{rp}] = [c_{pr}] = -c_{rp} \int_0^l N(x)^T N(x) dx$$

$$[c_{pp}]=(c_{rp}+c_p)\int_0^l N(x)^T N(x)dx$$

$$[c_{rr}]=c_{rp}\int_0^l N(x)^T N(x)dx+c_1 N(s)^T N(s)+2vm_e N(s)^T N'(s)$$

Stiffness Matrix component:

$$[k_{11}]=k_1$$

$$[k_{1r}]=-k_1 N(s)-c_1 vN'(s)$$

$$\{k_{r1}\}=k_1 N(s)^T$$

$$[k_{rp}]=[k_{pr}]=-k_{rp}\int_0^l N(x)^T N(x)dx$$

$$[k_{pp}]=E_p I_p\int_0^l N''(x)^T N''(x)dx+k_{rp}\int_0^l N(x)^T N(x)dx$$

$$[k_{rr}]=E_r I_r\int_0^l N''(x)^T N''(x)dx+k_{rp}\int_0^l N(x)^T N(x)dx$$
$$+k_1 N(s)^T N(s)+(c_1 v+m_e a)N(s)^T N'(s)+m_1 v^2 N(s)^T N''(s)$$

Force Vector components:

$$f_1=k_1 r(s)+c_1 vr'(s)$$

$$\{f_r\}=[(m_c+m_w)g-k_1-c_1 vr'(s)-m_e(ar'(s)+v^2 r''(s)]N(s)^T$$

$$\{f_p\}=\begin{Bmatrix}0\\0\\0\\0\end{Bmatrix}$$

$$\int_0^l N(x)^T N(x)dx=\begin{bmatrix}\frac{13l}{35} & \frac{11l^2}{210} & \frac{9l}{70} & \frac{-13l^2}{420}\\ \frac{11l^2}{210} & \frac{l^3}{105} & \frac{13l^2}{420} & -\frac{l^3}{140}\\ \frac{9l}{70} & \frac{13l^2}{420} & \frac{13l}{35} & -\frac{11l^2}{210}\\ \frac{-13l^2}{420} & -\frac{l^3}{140} & \frac{11l^2}{210} & \frac{l^3}{105}\end{bmatrix}$$

$$\int_0^l N''(x)^T N''(x)dx=\begin{bmatrix}\frac{12}{l^3} & \frac{6}{l^2} & \frac{-12}{l^3} & \frac{6}{l^2}\\ \frac{6}{l^2} & \frac{4}{l} & -\frac{6}{l^2} & \frac{2}{l}\\ \frac{-12}{l^3} & \frac{6}{l^2} & \frac{12}{l^3} & -\frac{6}{l^2}\\ \frac{6}{l^2} & \frac{2}{l} & \frac{-6}{l^2} & \frac{4}{l}\end{bmatrix}$$

### 2.3. Proposed Assembly Procedure

The proposed method consists assembling stiffness, mass and damping matrices, as well as the force vector of the system comprising rail and bridge, without any vehicle action.

#### 2.3.1. Bridge-Rail Elementary Interaction Matrix

Bridge and rail are each decomposed into 10 elements of equal length l, corresponding to 22 nodes. Each node has 2 DOFs Fig. 3, resulting in a total of 44 DOFs.

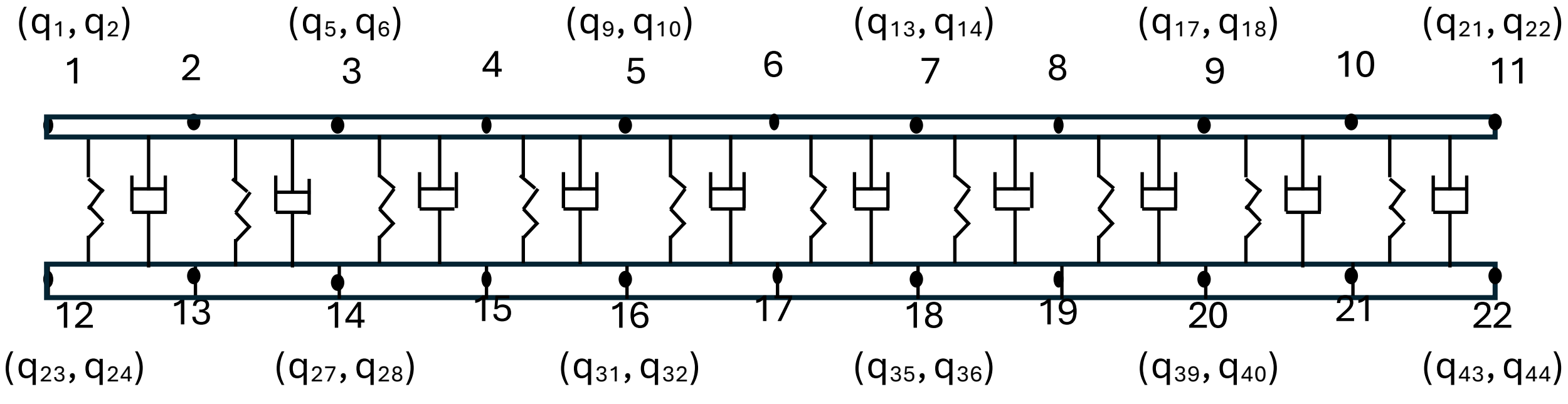

We note that:

**Figure 3.** Model of Rail- Bridge Interaction System without vehicle

The displacement vector of bridge-rail element system is 8 x 1, has an elementary matrix of 8 x 8 order.

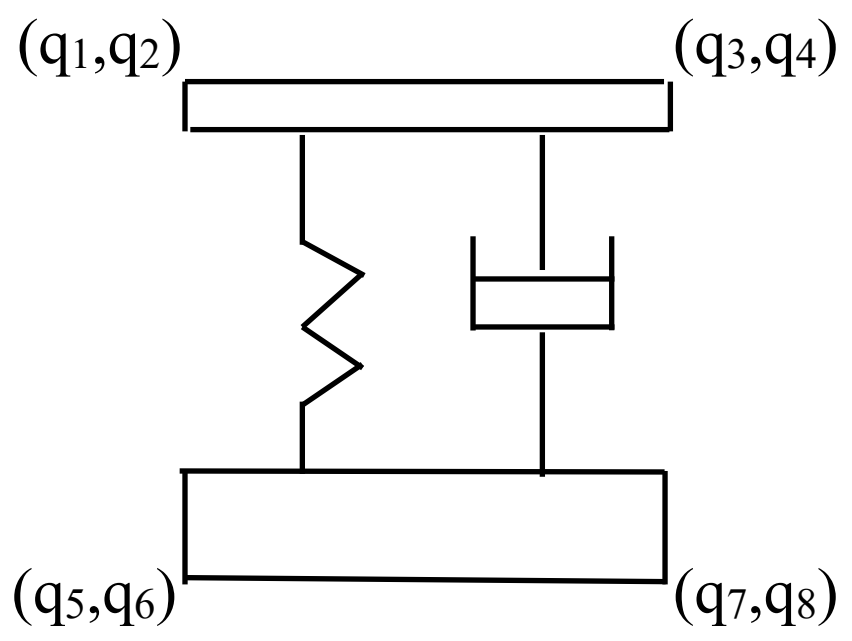

**Figure 4.** Model of Interaction Bridge-Track Element

The mass, stiffness and damping matrices as well as the force vector are defining interactions of the bridge-rail system for a single element:

- Elementary Mass Matrix $M_{rbe}$

$$M_{rbe}=\begin{bmatrix} [m_{rr}'] & [0] \\ [0] & [m_{pp}] \end{bmatrix}$$

$m_{rr}'$ and $m_{pp}$ blocks are defined as follow

$$[m_{rr}']=m_r\begin{bmatrix} \frac{13l}{35} & \frac{11l^2}{210} & \frac{9l}{70} & \frac{-13l^2}{420} \\ \frac{11l^2}{210} & \frac{l^3}{105} & \frac{13l^2}{420} & -\frac{l^3}{140} \\ \frac{9l}{70} & \frac{13l^2}{420} & \frac{13l}{35} & -\frac{11l^2}{210} \\ \frac{-13l^2}{420} & -\frac{l^3}{140} & \frac{11l^2}{210} & \frac{l^3}{105} \end{bmatrix}$$

$$[m_{pp}]=m_p\begin{bmatrix} \frac{13l}{35} & \frac{11l^2}{210} & \frac{9l}{70} & \frac{-13l^2}{420} \\ \frac{11l^2}{210} & \frac{l^3}{105} & \frac{13l^2}{420} & -\frac{l^3}{140} \\ \frac{9l}{70} & \frac{13l^2}{420} & \frac{13l}{35} & -\frac{11l^2}{210} \\ \frac{-13l^2}{420} & -\frac{l^3}{140} & \frac{11l^2}{210} & \frac{l^3}{105} \end{bmatrix}$$

- Elementary Stiffness Matrix $K_{rbe}$

$$K_{rbe}=\begin{bmatrix} [k_{rr}'] & [k_{rp}] \\ [k_{pr}] & [k_{pp}] \end{bmatrix}$$

$k_{rr}$ ‘, $k_{rp}$ and $k_{pp}$ blocks are defined as follow:

$$[k_{rr}']=E_rI_r\begin{bmatrix} \frac{12}{l^3} & \frac{6}{l^2} & \frac{-12}{l^3} & \frac{6}{l^2} \\ \frac{6}{l^2} & \frac{4}{l} & -\frac{6}{l^2} & \frac{2}{l} \\ \frac{-12}{l^3} & \frac{6}{l^2} & \frac{12}{l^3} & -\frac{6}{l^2} \\ \frac{6}{l^2} & \frac{2}{l} & \frac{-6}{l^2} & \frac{4}{l} \end{bmatrix}+k_{rp}\begin{bmatrix} \frac{13l}{35} & \frac{11l^2}{210} & \frac{9l}{70} & \frac{-13l^2}{420} \\ \frac{11l^2}{210} & \frac{l^3}{105} & \frac{13l^2}{420} & -\frac{l^3}{140} \\ \frac{9l}{70} & \frac{13l^2}{420} & \frac{13l}{35} & -\frac{11l^2}{210} \\ \frac{-13l^2}{420} & -\frac{l^3}{140} & \frac{11l^2}{210} & \frac{l^3}{105} \end{bmatrix}$$

$$[k_{rp}]=[k_{pr}]=-k_{rp}\begin{bmatrix} \frac{13l}{35} & \frac{11l^2}{210} & \frac{9l}{70} & \frac{-13l^2}{420} \\ \frac{11l^2}{210} & \frac{l^3}{105} & \frac{13l^2}{420} & -\frac{l^3}{140} \\ \frac{9l}{70} & \frac{13l^2}{420} & \frac{13l}{35} & -\frac{11l^2}{210} \\ \frac{-13l^2}{420} & -\frac{l^3}{140} & \frac{11l^2}{210} & \frac{l^3}{105} \end{bmatrix}$$

$$[k_{pp}]=E_pI_p\begin{bmatrix} \frac{12}{l^3} & \frac{6}{l^2} & \frac{-12}{l^3} & \frac{6}{l^2} \\ \frac{6}{l^2} & \frac{4}{l} & -\frac{6}{l^2} & \frac{2}{l} \\ \frac{-12}{l^3} & \frac{6}{l^2} & \frac{12}{l^3} & -\frac{6}{l^2} \\ \frac{6}{l^2} & \frac{2}{l} & \frac{-6}{l^2} & \frac{4}{l} \end{bmatrix}+k_{rp}\begin{bmatrix} \frac{13l}{35} & \frac{11l^2}{210} & \frac{9l}{70} & \frac{-13l^2}{420} \\ \frac{11l^2}{210} & \frac{l^3}{105} & \frac{13l^2}{420} & -\frac{l^3}{140} \\ \frac{9l}{70} & \frac{13l^2}{420} & \frac{13l}{35} & -\frac{11l^2}{210} \\ \frac{-13l^2}{420} & -\frac{l^3}{140} & \frac{11l^2}{210} & \frac{l^3}{105} \end{bmatrix}$$

- Elementary Damping Matrix $C_{rbe}$

$$[C_{rpe}]=\begin{bmatrix} [c_{rr}'] & [c_{rp}] \\ [c_{pr}] & [c_{pp}] \end{bmatrix}$$

$c_{rr}$ ‘, $c_{rp}$ and $c_{pp}$ blocks are defined as follow

$$[c_{rr}']=c_{rp}\begin{bmatrix} \frac{13l}{35} & \frac{11l^2}{210} & \frac{9l}{70} & \frac{-13l^2}{420} \\ \frac{11l^2}{210} & \frac{l^3}{105} & \frac{13l^2}{420} & -\frac{l^3}{140} \\ \frac{9l}{70} & \frac{13l^2}{420} & \frac{13l}{35} & -\frac{11l^2}{210} \\ \frac{-13l^2}{420} & -\frac{l^3}{140} & \frac{11l^2}{210} & \frac{l^3}{105} \end{bmatrix}$$

$$[c_{rp}]=[c_{pr}]=-c_{rp}\begin{bmatrix} \frac{13l}{35} & \frac{11l^2}{210} & \frac{9l}{70} & \frac{-13l^2}{420} \\ \frac{11l^2}{210} & \frac{l^3}{105} & \frac{13l^2}{420} & -\frac{l^3}{140} \\ \frac{9l}{70} & \frac{13l^2}{420} & \frac{13l}{35} & -\frac{11l^2}{210} \\ \frac{-13l^2}{420} & -\frac{l^3}{140} & \frac{11l^2}{210} & \frac{l^3}{105} \end{bmatrix}$$

$$[c_{pp}]=(c_{rp}+c_p)\begin{bmatrix} \frac{13l}{35} & \frac{11l^2}{210} & \frac{9l}{70} & \frac{-13l^2}{420} \\ \frac{11l^2}{210} & \frac{l^3}{105} & \frac{13l^2}{420} & -\frac{l^3}{140} \\ \frac{9l}{70} & \frac{13l^2}{420} & \frac{13l}{35} & -\frac{11l^2}{210} \\ \frac{-13l^2}{420} & -\frac{l^3}{140} & \frac{11l^2}{210} & \frac{l^3}{105} \end{bmatrix}$$

- Elementary Force Vector

$$\{f_p\}=\begin{Bmatrix} 0 \\ \vdots \\ 0 \end{Bmatrix}$$

2.3.2. Assembly of Elementary Matrix into Global Bridge-Rail System Matrix

The global mass, stiffness and damping matrix are 44 x 44 in size, and global displacement vector and force vector are 44 x 1 in size. Connectivity vectors for each element are defined as follows:

**Table 1.** Connectivity vectors for each element

| Element | Connectivity vectors | |
|---|---|---|
| 1 | VC1 | [1 2 3 4 23 24 25 26] |
| 2 | VC2 | [3 4 5 6 25 26 27 28] |
| 3 | VC3 | [5 6 7 8 27 28 29 30] |
| 4 | VC4 | [7 8 9 10 29 30 31 32] |
| 5 | VC5 | [9 10 11 12 31 32 33] |
| 6 | VC6 | [11 12 13 14 33 34 35 ] |
| 7 | VC7 | [13 14 15 16 35 36 37 ] |
| 8 | VC8 | [15 16 17 18 37 38 39 ] |
| 9 | VC9 | [17 18 19 20 39 40 41] |
| 10 | VC10 | [19 20 21 22 41 42 43] |

The assembly function connecting the two matrices is defined as follows: g being the elementary matrix to be assembled and G the overall matrix of the system.

For all (i,j) ∈[1 2 3 4 5 6 7 8 9 10 ]2

$$\begin{cases} G(VCp(i),VCp(j)) = g(i,j) & p=1 \\ G(VCp(i),VCp(j)) = G(VCp(i),VCp(j)) + g(i,j) & p\neq 1 \end{cases}$$

VCp is element p connectivity vector.

The global matrices of the bridge-rail system are obtained by assembling the elementary matrices $[C_{rpe}]$, $[k_{rpe}]$, $[M_{rpe}]$ using the function defined above. Therefore, we find respectively matrices of the bridge-rail system without any vehicle action $[C_{rpg}]$, $[k_{rpg}]$ and $[M_{rpg}]$.

2.3.3. Matrix of "Vehicle-Rail" Interaction

- Mass Matrix

Consider $[M_{rv}]$ the mass matrix induced by the moving mass on the rail, with size 44 x 44, s denotes the local coordinate measured from the left node of a beam element.

$$[M_{rv}]=m_w\times N(s)^T\times N(s)$$

$$N(s)=[0 \quad 0 \quad 0 \quad N_A(s)\ N_B(s)\ N_c(s)\ N_D(s) \quad 0 \quad 0 \quad 0]$$

N(s) is of size 44 x 1 and all its elements are null except those which correspond to rail element DOF on which the wheel operates.

- Stiffness Matrix

Consider$[K_{rv}]$ , $[k_{1r}]$ and $\{k_{r1}\}$ are the stiffness matrices induced by wheel displacement on rail, they are respectively of order 44×44, 1 × 44 and 44×1 with:

$$[K_{rv}]=k_1N(s)^TN(s)+(c_1v+m_ea)N(s)^TN'(s(s)) +m_1v^2N(s)^TN''(s)$$

$$[k_{1r}]=-k_1N(s)-c_1vN'(s)$$

$$\{k_{r1}\}=k_1N(s)^T$$

- Damping Matrix

The matrices $[C_{rv}]$ , $[c_{1r}]$ and $\{c_{r1}\}$ are of orders 44 × 44 ,1×44 et 44×1 respectively and represent overall damping matrices induced by the wheel on rail.

$$[C_{rv}]=c_1N(\xi)^TN'(\xi)+2vm_eN(\xi)^TN'(\xi)$$

$$[c_{1r}]=-c_1N(s)$$

$$\{c_{r1}\}=-c_1N(s)$$

- Force Vector

$\{f_{rv}\}$ represent the load force vector induced by the wheel acting on the rail, it is of order 44×1:

$$\{f_{rv}\}=[(m_c+m_w)g-k_1-c_1vr'-m_e(ar'+v^2r'']N((s))^T$$

2.3.4. Construction of Global Matrices for the Vehicle-Rail-Axle System

The global matrices of the system composed of 3 elements: vehicle, rail and bridge are written as follows:

- Mass Matrix:

$$M_g=\begin{bmatrix} [m_{11}] & [0] \\ \{0\} & [m_{rbe}]+[m_{rv}] \end{bmatrix}$$

- Damping Matrix:

$$C_g=\begin{bmatrix} [c_{11}] & [c_{1r}] \\ \{c_{r1}\} & [c_{rbe}]+[c_{rv}] \end{bmatrix}$$

- Stiffness Matrix:

$$K_g=\begin{bmatrix} [k_{11}] & [k_{1r}] \\ \{k_{r1}\} & [k_{rbe}]+[k_{rv}] \end{bmatrix}$$

- Force Vector:

$$F_g=\begin{Bmatrix} f_1 \\ \{f_{rv}\}+\{f_p\} \end{Bmatrix}$$

The dynamic response vectors are respectively:

- Displacement:

$$Q_g=\begin{Bmatrix} q_1 \\ \{q_r\} \\ \{q_p\} \end{Bmatrix}$$

- Velocity:

$$\dot{Q}_g=\begin{Bmatrix} \dot{q}_1 \\ \{\dot{q}_r\} \\ \{\dot{q}_p\} \end{Bmatrix}$$

- Acceleration:

$$\ddot{Q}_g=\begin{Bmatrix} \ddot{q}_1 \\ \{\ddot{q}_r\} \\ \{\ddot{q}_p\} \end{Bmatrix}$$

for It should be noted that:

- $\{\ddot{q}_r\}, \{\dot{q}_r\}$ et $\{q_r\}$ are of order $22\times 1$. (The rail's degrees of freedom).
- $\{\ddot{q}_p\}, \{\dot{q}_p\}$ et $\{q_p\}$ are of order $22 \times 1$. (The bridge's degrees of freedom)
- $\{f_{rv}\}$ is of order $44 \times 1$

The assembly process used allows genering the global motion equation for the entire vehicle-rail-bridge interaction system follows:

$$M_g\ddot{Q}_g+K_gQ_g+C_g\dot{Q}_g=F_g \tag{16}$$

# 3. Dynamic Analysis Using Newmark Beta Method

Several numerical integration methods in the time domain [20] have been reported to solve these types of equation. Three main requirements must be addressed in a numerical solution procedure, namely:

- Convergence: the numerical solution approaches the exact solution as the time step decreases.

- Accuracy: the numerical solution presents result close enough to exact solution.

- Stability: the solution must be stable even in presence of errors.

The Newmark beta method remains one of the most popular explicit methods used in structures dynamic analysis. It is applicable to linear differential systems with time-dependent mass, stiffness and damping matrices [21].

In addition to $\beta$ parameter, it contains a $\delta$ parameter with the value ½. The Newmark equation can be written in incremental quantities for a constant time step Δt as follows:

$$\Delta u_i=\dot{u}_i\Delta t+\frac{1}{2}\ddot{u}_i\Delta t^2+\beta\Delta\ddot{u}_i\Delta t^2 \tag{17}$$

$$\Delta\dot{u}_i=\ddot{u}_i\Delta t+\delta\Delta\ddot{u}_i\Delta t \tag{18}$$

$$\Delta u=u(t+\Delta t)-u(t) \tag{19}$$

$$\Delta\dot{u}=\dot{u}(t+\Delta t)-\dot{u}(t) \tag{20}$$

$$\Delta\ddot{u}=(\Delta u-\dot{u}_t\Delta t-\frac{1}{2}\ddot{u}_t\Delta t^2)/\beta\Delta t^2 \tag{21}$$

Substituting "(21) "into "(18)"

$$\Delta\dot{u}=\frac{1}{2\beta\Delta t}\Delta u-\frac{1}{2\beta}\dot{u}_t+\left(1-\frac{1}{4\beta}\right)\ddot{u}_t\Delta t^2) \tag{22}$$

Considering the equation of motion at time t and time t+Δt.

$$M\ddot{u}_t+C\dot{u}_t+Ku_t=F_t \tag{23}$$

$$M\ddot{u}_{t+\Delta t}+C\dot{u}_{t+\Delta t}+Ku_{t+\Delta t}=F_{t+\Delta t} \tag{24}$$

Subtracting "(22)" from "(23)":

$$M\Delta\ddot{u}+C\Delta\dot{u}+K\Delta u=\Delta F \tag{25}$$

With:

$$\Delta\ddot{u}=\ddot{u}_{t+\Delta t}-\ddot{u}_t \tag{26}$$

$$\Delta\dot{u}=\dot{u}_{t+\Delta t}-\dot{u}_t \tag{27}$$

$$\Delta u=u_{t+\Delta t}-u_t \tag{28}$$

The values of M, C, and K in " (25) "are calculated at time t, assuming that remain constant during Δt.

Substituting "(21)" and "(22)" into "(25)" gives "(30)", which calculates the displacement Δu:

$$\widehat{K_t}\Delta u=\widehat{\Delta F} \tag{29}$$

Where the effective stiffness $\widehat{K_t}$and incremental force $\widehat{\Delta F}$ are provided respectively

$$\widehat{K_t}=K_t+\frac{1}{\beta\Delta t^2}M_t+\frac{1}{2\beta\Delta t}C_t \tag{30}$$

$$\widehat{\Delta F}=\Delta F+\frac{1}{\beta\Delta t}M_t\dot{u}_t+\frac{1}{2\beta}C_t\dot{u}_t+\frac{1}{2\beta}M_t\ddot{u}_t-C_t\Delta t\left(1-\frac{1}{4\beta}\right)\ddot{u}_t \tag{31}$$

The value of the β parameter should be between 1/6 and 1/2.

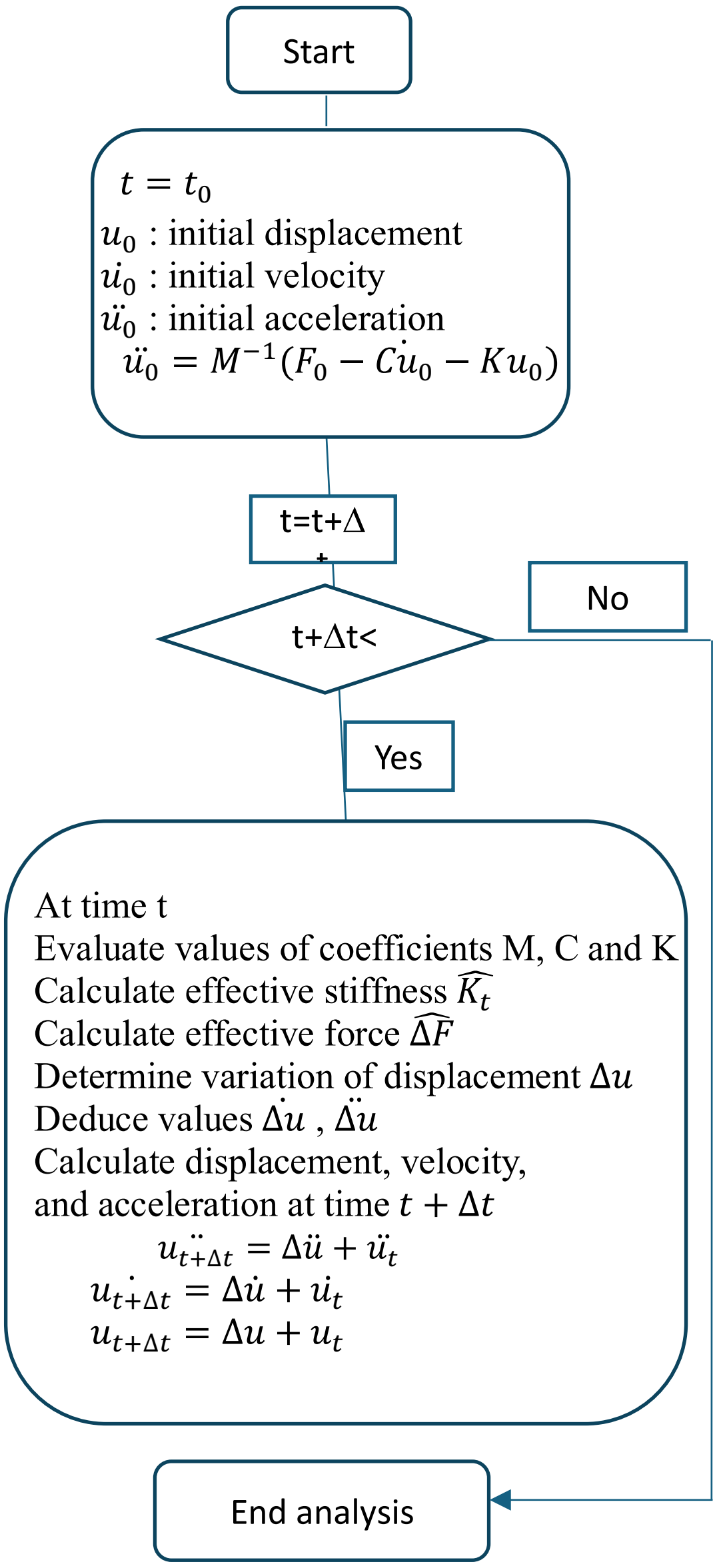

**Figure 5.** Flowchart for process of equation resolution

# 4. Verification of Proposed Procedure

The assembly procedure used to obtain the motion equations of the vehicle-rail-bridge interaction system and the associated computer program are verified through the study of an example in this section.

The differential equations were solved by the Newmark-β method with a time step Δt = 0.005. The computational program was developed on "MATLAB © R2021b, The MathWorks, Inc.".

The example studied to validate the proposed procedure consists of a simply supported Bernoulli beam of length L = 30 m, travelling with a single-axle vehicle at constant speed v = 27.78 m/s. The damping ratio of the bridge

ζ = 0 and the surface of the beam is assumed to be smooth. Vehicle characteristics, bridge and rail are summarized in Table 2 [1].

**Table 2.** System vehicle – rail-bridge parameters

| Notation | Value | Unit |
|---|---|---|
| Vehicle | | |
| $m_1$ | 5750 | Kg |
| $m_e$ | 0 | Kg |
| $C_v$ | 0 | Ns/m |
| $K_v$ | $1.595\ 10^6$ | N/m |
| Rail | | |
| $m_r$ | $10^{-7}$ | Kg/m |
| $I_r$ | $10^{-10}$ | $m^4$ |
| $E_r$ | $2.06\ 10^{11}$ | Pa |
| $K_{rb}$ | $10^{13}$ | N/m |
| $C_{rb}$ | 0 | Ns/m |
| Bridge | | |
| $m_b$ | $2.303\ 10^3$ | Kg/m |
| $I_b$ | 2.90 | $m^4$ |
| $E_b$ | $2.87\ 10^9$ | Pa |
| $C_b$ | 0 | Ns/m |

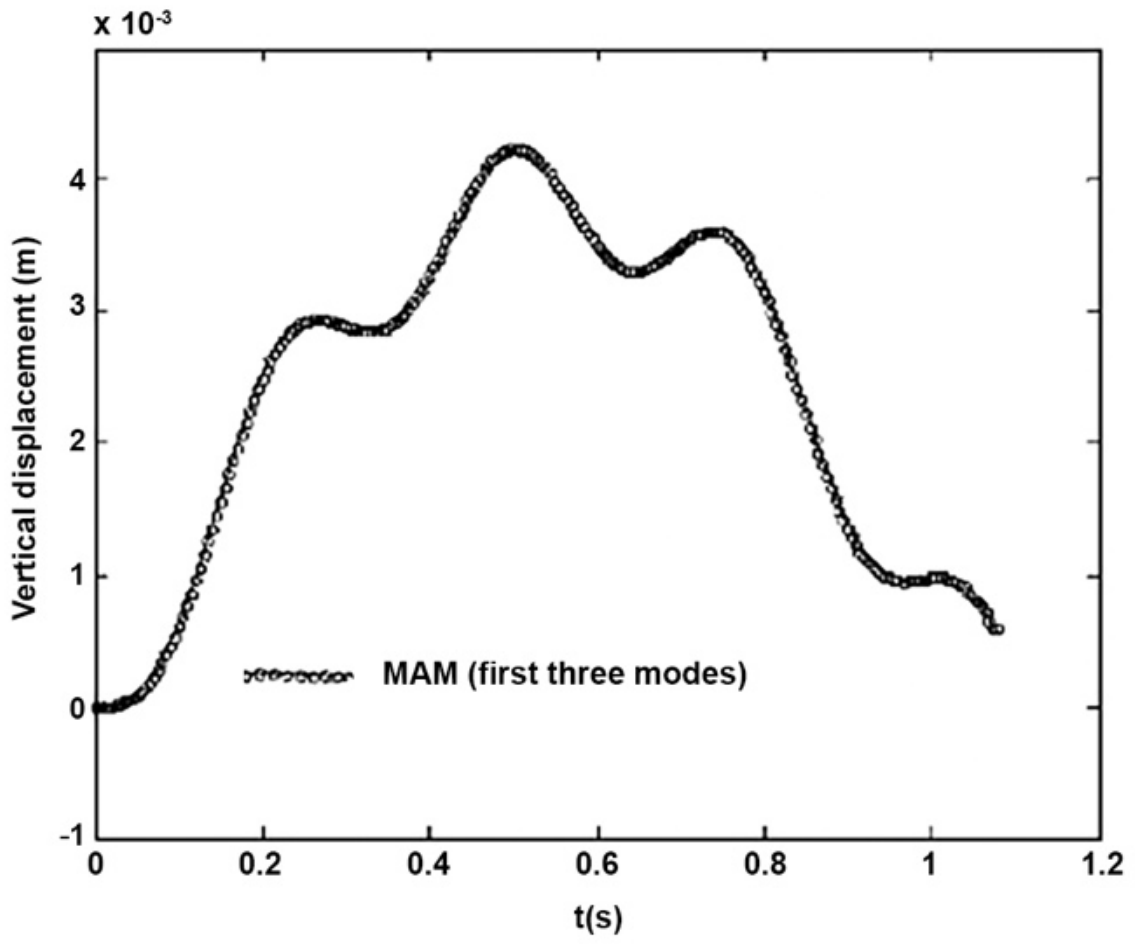

**Figure 6.** Time history of vertical displacement of midpoint of bridge computed by modal analysis method [1]

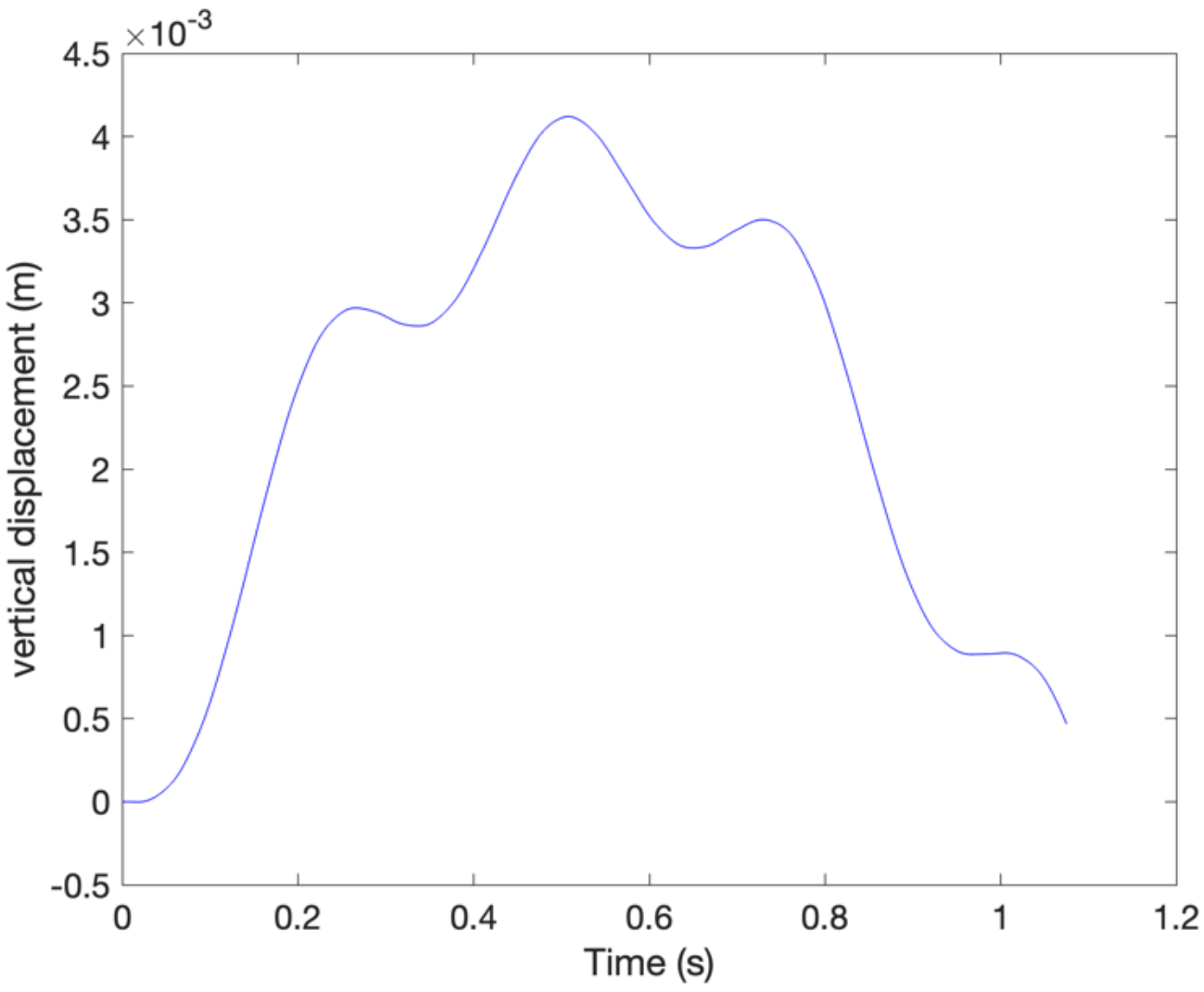

**Figure 7.** Time history of vertical displacement of midpoint of bridge calculated by current method

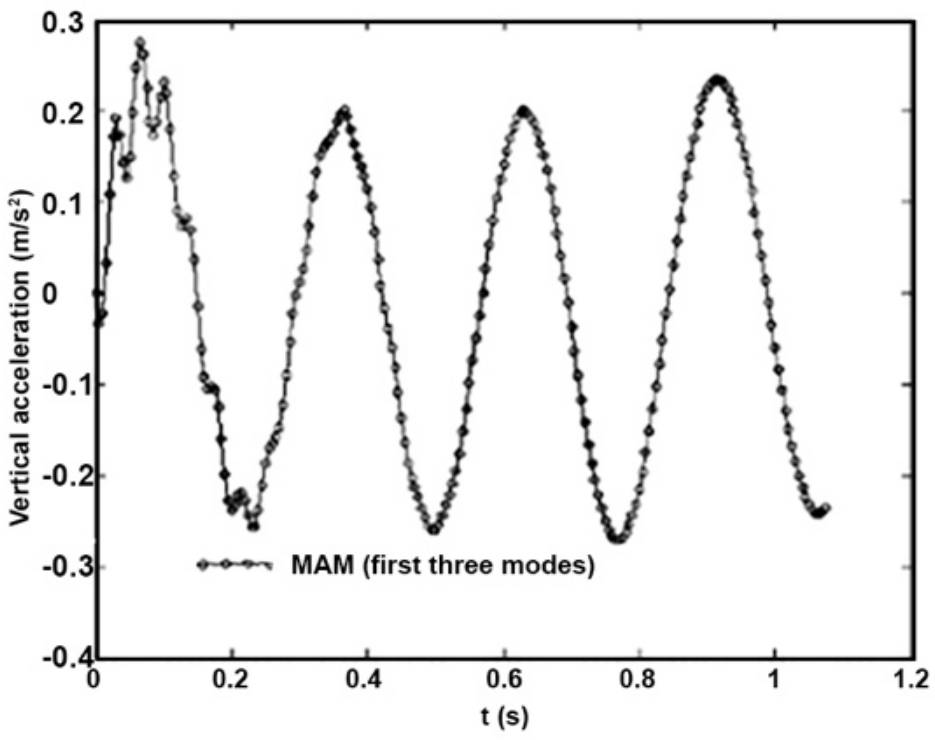

**Figure 8.** Time history of vertical acceleration of midpoint of bridge computed by analysis method [1]

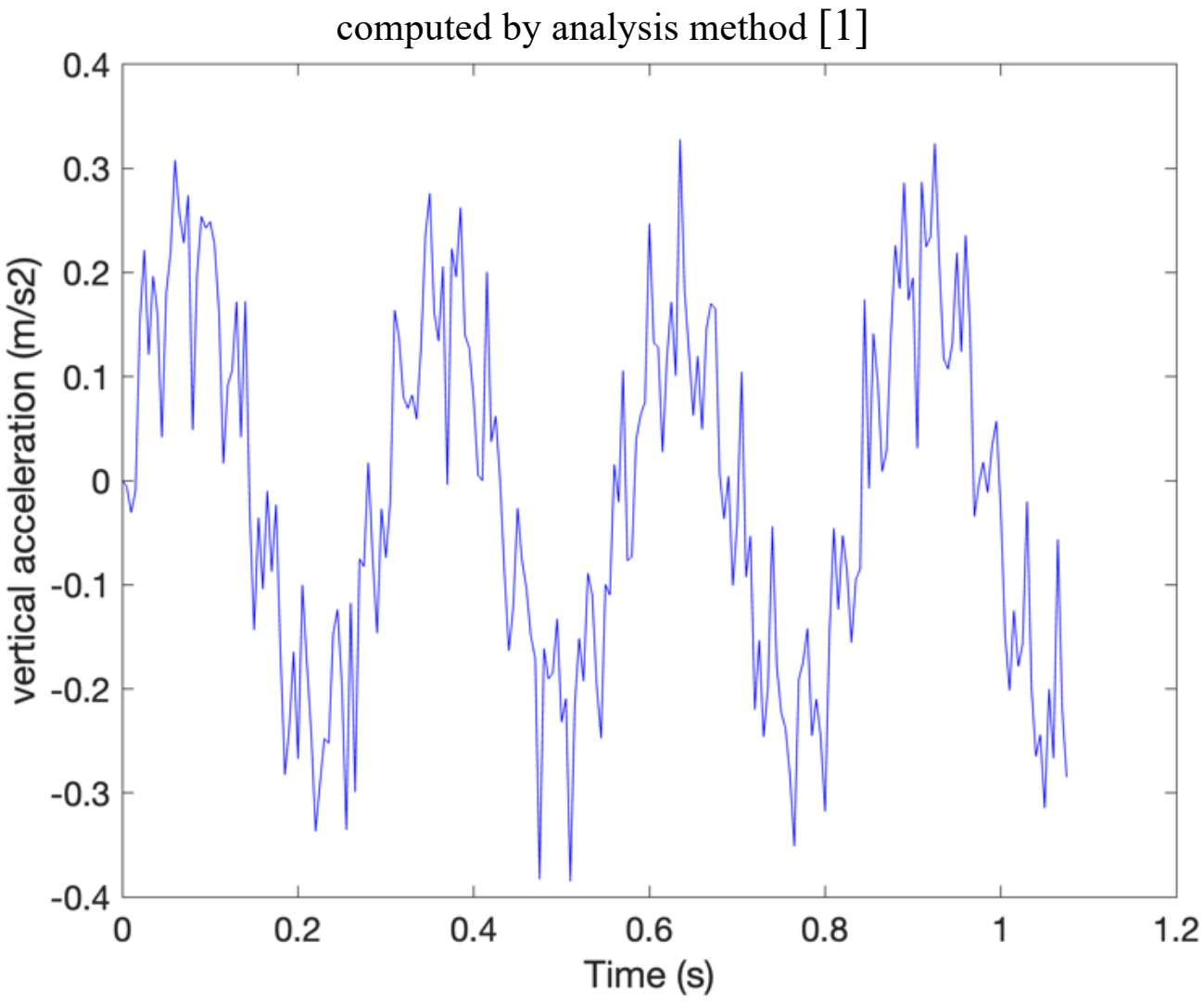

**Figure 9.** Time history of vertical acceleration of midpoint of bridge calculated by current method

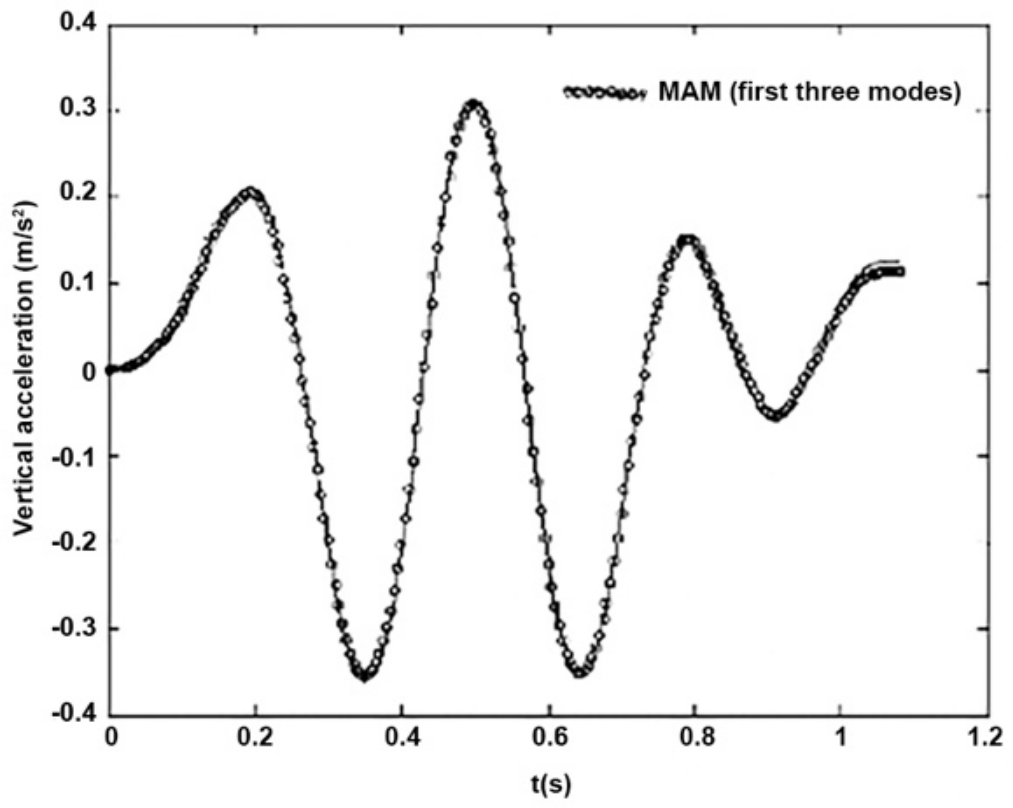

**Figure 10.** Time history of vertical acceleration of carbody computed by modal analysis method [1]

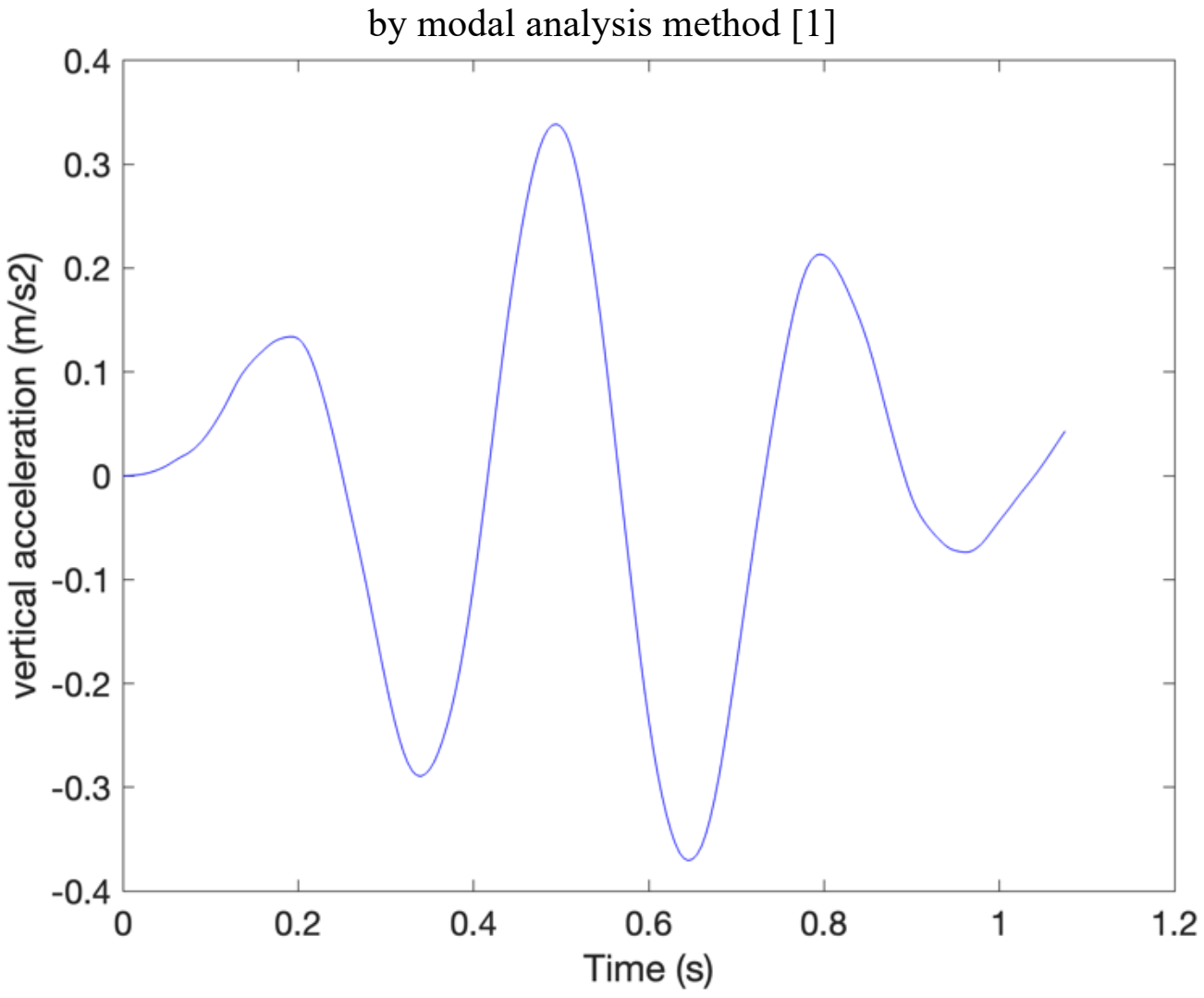

**Figure 11.** Time history of vertical acceleration of carbody calculated by current method

# 5. Conclusion and Remarks

Based on the energy approach and the finite element method, the motion equations of a vehicle-rail-bridge system element have been established. The model presented divides the bridge and rail into 10 elements of equal length l. The vehicle is modeled by a suspended mass moving along the 10 elements at a constant speed v. The elements to be assembled are not identical which increases the complexity of their assembly. In fact, there are two types of elements to be assembled: an element with a vehicle and 9 others without a vehicle.

This study elaborates a procedure for assembling elementary equations into the system global equations under study. These equations are then solved using the explicit Newmark beta numerical method to obtain the dynamic responses of the overall vehicle-rail-bridge system.

The procedure used was validated by comparing its results with those obtained by the method (MAM) presented in article [1]. The method presented is however restricted to bridges of limited length, since it only defines assembly of 10 elements. A more advanced technique for assembling n elements of the elementary system is currently under investigation, based on a vehicle model with 4 wheels and 10 degrees of freedom.